\documentclass{amsart}

\usepackage{amsmath}
\usepackage{amssymb}
\usepackage{graphicx}
\usepackage{mathrsfs}
\usepackage{tikz-cd}
\usepackage{xcolor}
\usepackage{tikz}
\usepackage[utf8]{inputenc}
\usepackage[english]{babel}
\usepackage{amsthm}
\usepackage{comment}

\theoremstyle{plain}
\newtheorem{theorem}[subsection]{Theorem}
\newtheorem{lemma}[subsection]{Lemma}
\newtheorem{corollary}[subsection]{Corollary}
\newtheorem{proposition}[subsection]{Proposition}
\theoremstyle{definition}

\newtheorem{definition}[subsection]{Definition}
\newtheorem{remark}[subsection]{Remark}

\usetikzlibrary{%
  matrix,%
  calc,%
  arrows%
}

\graphicspath{ {Images/} }

\newcommand{\C}{\mathbb{C}}
\newcommand{\R}{\mathbb{R}}
\newcommand{\Z}{\mathbb{Z}}
\newcommand{\Q}{\mathbb{Q}}
\newcommand{\N}{\mathbb{N}}

\newcommand{\tensor}{\otimes}
\DeclareMathOperator{\supp}{supp}
\DeclareMathOperator{\Tr}{Tr}
\DeclareMathOperator{\id}{id}

\DeclareMathOperator{\re}{Re}

\DeclareMathOperator{\ch}{ch}
\DeclareMathOperator{\diag}{diag}
\DeclareMathOperator{\vecc}{Vec}
\DeclareMathOperator{\End}{End}

\usepackage{esint}

\usepackage[backref=page, colorlinks]{hyperref}

\title[Vanishing of Topological Invariants]{Vanishing of Topological Invariants For Unnormalized Schatten $p$ multiplicative Maps}

\author{Forrest Glebe}

\begin{document}

\begin{abstract}
A result of Dadarlat shows that nonzero even rational cohomology obstructs the matricial stability of many discrete groups. In the author's previous work, 2-cohomology is used to argue that certain groups are not stable in unnormalized Schatten $p$-norms for $p > 1$, although 2-cohomology is known not to obstruct stability in the unnormalized 2-norm in general. The main result of this paper demonstrates that we should not expect $2k$ cohomology to obstruct in the unnormalized Schatten $p$-norm for $p\le k$, because the invariant in Dadarlat's argument vanishes for maps that are asymptotically multiplicative in the Schatten $p$-norm.
\end{abstract}

\maketitle

\section{Introduction}

A countable discrete group $\Gamma$ is said to be {\em stable in the unnormalized Schatten $p$-norm} if every function from the group that is ``almost multiplicative'' in the point $p$-norm topology is ``close'' to a genuine unitary representation in the same topology; see Definition~\ref{stability} for a formal description. Of particular interest is the $p=2$ case called {\em Frobenius stability}. Stability in the operator norm (which we will also refer to as the $\infty$-norm) is called {\em matricial stability}.

Frobenius stability was introduced by de Chiffre, Glebsky, Lubotzky, and Thom in~\cite{nonaproximable}. Stability of a finitely presented group is equivalent to a notion of stability of the presentation of that group; this notion was shown to be independent of the presentation by Arzhantseva and Păunescu in~\cite{almostperm}.

A result of Dadarlat shows that if a $H^{2k}(\Gamma;\Q)\ne0$ for some $k>0$ and $\Gamma$ satisfies some additional mild but technical conditions, then $\Gamma$ is not matrically stable~\cite{obs}. Previous work by the author uses 2-cohomology to argue non-stability of certain groups in the unnormalized Schatten $p$-norm for $p>1$~\cite{frob}~\cite{cup}, but in general, neither 2-cohomology nor higher even cohomology obstructs Frobenius stability~\cite{lattice}. This raises the question of whether or not higher cohomology can be used to argue for the non-stability of any groups in unnormalized Schatten $p$-norms. To make this question more precise, we describe the invariants that come up in Dadarlat's proof.

Suppose $\rho_n:\Gamma\rightarrow U(m_n)$ is an operator norm asymptotically multiplicative sequence of maps, and suppose that $Y\subseteq B\Gamma$ is compact. For large enough $n$, following Phillips and Stone \cite{phillips-stone1}\cite{phillips-stone2}, we may associate a vector bundle $E_{\rho_n}^Y$ over $Y$ (see Definition~\ref{E_rho}). Then the rational Chern Character $\ch_k(E_{\rho_n}^Y)\in H^{2k}(Y;\Q)$ can be used to show that the asymptotic homomorphism is far from a genuine representation; if there exists a $k>0$ so that for all $n$, the $\ch_k(E_{\rho_n}^Y)\ne0$, then $\rho_n$ cannot be approximated by genuine representations. This is equivalent to the invariant used in \cite{obs}. The the ``winding number argument'' of Kazhdan~\cite{KvoiculescuM}, and Exel and Loring~\cite{ELVoiculescuM}, is equivalent to the pairing of $\ch_1(E_{\rho_n}^Y)$ with a homology class of $\Gamma$. The invariant used in~\cite{frob} and~\cite{cup} is, in turn, equivalent to the ``winding number argument.'' These equivalences are shown in~\cite{partB} and~\cite{partA}, respectively. We aim to show the following:
\begin{theorem}\label{main} 
Suppose that $\Gamma$ is a countable discrete group, $\rho_n:\Gamma\rightarrow U(m_n)$ is a sequence of functions satisfying
$$||\rho_n(gh)-\rho_n(g)\rho_n(h)||_p\rightarrow0\quad\forall g,h\in\Gamma$$
where $||\cdot||_p$ is the unnormalized Schatten $p$-norm, and $Y\subseteq B\Gamma$ is compact. If $k\ge p$ then for sufficiently large $n$, $\ch_k(E_{\rho_n}^Y)=0$.
\end{theorem}

As a result, we should not expect $2k$ cohomology to obstruct stability in the unnormalized Schatten $p$-norm when $k\ge p$, at least not using the same obstruction in~\cite{obs}. In particular, we should not expect higher even cohomology to obstruct Frobenius stability. 

The bound $k\ge p$ is sharp by a recent example in joint work of Dadarlat and the author~\cite[Remark 8.6]{partB}; there is a sequence of functions $\rho_n$ from $\Z^{2k}$ to unitary matrices so that $\ch_k(E_{\rho_n}^{B\Z^{2k}})\ne0$, but $\rho_n$ is asymptotically multiplicative in the unnormalized Schatten $p$-norm for all $p>k$. Indeed, one can pick the maps so that $\ch_j(E_{\rho_n}^{B\Z^{2k}})=0$ for all $0<k<j$. Theorem~\ref{main} would not hold if Chern character were replaced with the Chern class (Remark~\ref {chernclass}).

For readers interested in operator norm stability Theorem~\ref{main} implies the following Corollary:

\begin{corollary}\label{dimbound}
Let $\Gamma$ be a countable discrete group and $Y\subseteq B\Gamma$ be compact, let $k,N\in\N_{>0}$. Then there exists a finite subset $S\subseteq\Gamma$ and, $\varepsilon>0$ so that for all $n\in\N_{>0}$ if $\rho:\Gamma\rightarrow U(Nn^k)$ and $||\rho(g)\rho(h)-\rho(gh)||_\infty<\frac{\varepsilon}{n}$ for all $g,h\in S$ it follows that $\ch_k(E_\rho^Y)=0$.
\end{corollary}
In other words, if $\rho_n$ is an asymptotic homomorphism with nonzero $k$th Chern character, either the dimension grows more quickly than $Nn^k$ or the failure to be multiplicative goes to zero at most as quickly as $\frac{1}{n}$.

Another application is to show that there are asymptotically multiplicative maps in operator norm that cannot be ``improved'' to a map that is asymptotically multiplicative in the unnormalized Schatten $p$-norms. 

In our notation, as in~\cite{obs}, we define a group to be {\em quasi-diagonal} if it admits a faithful quasi-diagonal unitary representation. All residually amenable groups are known to be quasi-diagonal \cite[Corollary C]{nuclearQD}\cite[Remark 3.5]{obs}. 
The definition of a $\gamma$-element is technical; see~\cite{KKNovikov}. 
We note that groups with Haagerup's property~\cite{embeddable_gamma}\cite{haagerup} (in particular amenable groups), linear groups~\cite{embeddable_gamma}\cite{linear_gamma}, and hyperbolic groups~\cite{embeddable_gamma}\cite{brownozawa} all have a $\gamma$-element. 
See,~\cite{buildings} and~\cite{bolic} for more examples.

\begin{corollary}\label{can'timprove}
If $\Gamma$ is a countable quasidiagonal group with a gamma element so that $H^{2k}(\Gamma,\R)\ne0$, and $1\le p\le k$, there is a sequence of functions $\psi_n:\Gamma\rightarrow U(m_n)$, satisfying
$$||\psi_n(g)\psi_n(h)-\psi_n(gh)||_\infty\rightarrow0$$
for all $g,h\in\Gamma$, but there is no sequence of functions $\rho_n:\Gamma\rightarrow U(m_n)$ so that
$$||\psi_n(g)-\rho_n(g)||_\infty\rightarrow0$$
and
$$||\rho_n(g)\rho_n(h)-\rho_n(gh)||_p\rightarrow0$$
for all $g,h\in\Gamma$.
\end{corollary}

Finally, we show the following $p$-norm to operator norm stability condition for certain groups with vanishing low-dimensional even cohomology.

\begin{corollary}\label{cor3}
Let $\Gamma$ be a countable residually-finite torsion-free amenable group and let $1\le p<\infty$. Suppose that for $0<k<p$, $H^{2k}(\Gamma;\Q)=\{0\}$. If $\rho_n:\Gamma\rightarrow U(m_n)$ is a sequence of functions so that
$$||\rho_n(gh)-\rho_n(g)\rho_n(h)||_p\rightarrow0\quad\forall g,h\in\Gamma.$$
Then there exist sequences of representations $\pi_n$, $\psi_n$ and integers $r_n$ so that
$$||\rho_n(g)^{\oplus r_n}\oplus\psi_n(g)-\pi_n(g)||_\infty\rightarrow0\quad\forall g\in\Gamma.$$
\end{corollary}

The paper is organized as follows. In Section~2, we cover relevant preliminaries. In Section~3, we introduce notation for matrices of differential forms and prove some basic relevant facts. In Section~4, we prove the main result in the case that $Y$ can be given a smooth structure (Theorem~\ref{main2}). In Section~5, we deduce the main result in general and prove the other Corollaries.

There is a well-known correspondence between operator norm almost representations of $\Gamma$ and operator norm almost-flat bundles over $B\Gamma$, or compact subsets of $B\Gamma$~\cite{CGM}\cite{phillips-stone1}\cite{phillips-stone2}\cite{manuilovmishchenko}\cite{mishchenkoteleman}\cite{quasi-rep}\cite{mariusjosealmostflat}. Here, an almost-flat bundle can be described as either a bundle with almost-constant cocycles or with a small curvature form~\cite{mishchenkoteleman}. The key idea behind the proof is to find an appropriate definition of the $p$-norm of the curvature form (Definition~\ref{norms}), then show that if a map is almost multiplicative in the $p$-norm, the curvature form is small in the $p$-norm. From this, we use H\"older's inequality and a well-known formula for Chern characters (Proposition~\ref{chernFormula}) to show that for $k>p$ the $k$th Chern character is the de Rham cohomology class of a form that is uniformly almost zero, and thus the class must be zero.

\section{Preliminaries}

\begin{definition}\label{stability}
Let $\Gamma$ be a countable discrete group and let $1\le p\le\infty$. Let $||\cdot||_p$ denote the unnormalized Schatten $p$-norm, $||M||_p=(\Tr((M^*M)^{p/2}))^{1/p}$ for $p<\infty$ and operator norm for $p=\infty$. Let $U(k_n)$ be the $k_n\times k_n$ complex unitary group.  A {\em asymptotic homomorphism in the unnormalized Schatten $p$-norm} is a sequence of functions $\rho_n:\Gamma\rightarrow U(k_n)$ so that for all $g,h\in\Gamma$ we have 
$$||\rho_n(gh)-\rho_n(g)\rho_n(h)||_p\rightarrow0,$$ We say $\Gamma$ is {\em stable in the unnormalized Schatten $p$-norm} if every for every asymptotic homomorphism in the unnormalized Schatten $p$-norm, $\rho_n:\Gamma\rightarrow U(k_n)$ there a sequence of unitary representations $\psi_n:\Gamma\rightarrow U(k_n)$ so that
$$||\rho_n(g)-\psi_n(g)||_p\rightarrow0$$
for all $g\in\Gamma$.
\end{definition}

\begin{lemma}\label{calculation}
Let $z$ be any complex number.
\begin{enumerate}
    \item Let $S=\{x\in\C:\re(x)\ge\frac12\}$. Then
$$|\chi_S(z)-z|\le 2|z^2-z|.$$
    \item Let $R=\{x\in\C:\re(x)>0\}$. Then $|z-(2\chi_R(z)-1)|\le|z^2-1|.$
\end{enumerate}
\end{lemma}
\begin{proof}
(1) We compute,
\begin{align*}
|z-\chi_{S}(z)|&=\min\{|z|,|z-1|\}\\
&\le 2\min\{|z|,|z-1|\}\cdot\max\{|z|,|z-1|\}\\
&=2|z|\cdot|z-1|\\
&=2|z^2-z|
\end{align*}
The inequality in the second line holds since at least one $|z|$ and $|z-1|$ will be at least $\frac12$ for all $z$.

(2) Note that $|z-(2\chi_R(z)-1)|=\min\{|z-1|,|z+1|\}$. Then the proof follows the proof of (1).
\end{proof}

\begin{definition}
We say that a function $\rho_n:\Gamma\rightarrow U(k_n)$ is {\em normalized} if $\rho_n(1)=1$ and $\rho_n(g^{-1})=\rho_n(g)^*$ for all $g\in\Gamma$.
\end{definition}

\begin{proposition}\label{balanced}
Suppose that $\rho_n:\Gamma\rightarrow U(k_n)$ is asymptotically multiplicative in the unnormalized Schatten $p$-norm. Then there exists a sequence of normalized functions $\psi_n:\Gamma\rightarrow U(k_n)$ so that $||\psi_n(g)-\rho_n(g)||_p\rightarrow0$ for all $g\in\Gamma$.
\end{proposition}
\begin{proof}
Using the axiom of choice we may partition $\Gamma$ as follows,
$$\Gamma=\{e\}\sqcup\Gamma_1\sqcup\Gamma_2\sqcup\Gamma_3,$$
where $\Gamma_1$ consists of all 2-torsion elements and a non-2-torsion non-identity element $g\in\Gamma_2$ if and only if $g^{-1}\in\Gamma_3$. Let $R=\{z\in\C:\re(z)>0\}$ and let $\chi_{R}$ be the characteristic function of $R$. Then let
$$\psi_n(g)=\begin{cases}
1&g=e\\2\chi_{R}(\rho_n(g))-1&g\in\Gamma_1\\\rho_n(g)&g\in\Gamma_2\\\rho_n(g^{-1})^*&g\in\Gamma_3.
\end{cases}$$
Clearly, $\psi_n$ is normalized. First, we see that
$$||\rho_n(e)-1||_p=||\rho_n(e)^2-\rho_n(e)||_p$$
which goes to zero because $\rho_n$ is asymptotically multiplicative. Next let $g\in \Gamma_1$. Pick a basis of $\C^{k_n}$ so that $\rho_n(g)$ is diagonal. Then $\rho_n(g)=\diag(\lambda_1,\ldots,\lambda_{k_n})$. Using Lemma~\ref{calculation} we see that
\begin{align*}
||\rho_n(g)-\psi_n(g)||_p&=||\langle\lambda_1-(2\chi_R(\lambda_1)-1),\ldots,\lambda_{k_n}-(2\chi_R(\lambda_{k_n})-1)\rangle||_p\\
&\le||\langle\lambda_1^2-1,\ldots\lambda_{k_n}^2-1\rangle||_p\\
&=||\rho_n(g)^2-1||_p\\
&\le||\rho_n(g)^2-\rho_n(e)||+||\rho_n(e)-1||_p.
\end{align*}
The first term goes to zero by asymptotic multiplicativity and we have already shown that the second term goes to 0. Finally, let $g\in\Gamma_3$. We compute
\begin{align*}
||\rho_n(g)-\rho_n(g^{-1})^*||_p&=||\rho_n(g)\rho_n(g^{-1})-1||_p.
\end{align*}
The term on the right goes to zero for the same reasoning expressed above. 
\end{proof}

The proof above applies to operator norm and normalized Schatten $p$-norms as well.

\begin{definition}
A tool we will use frequently is the holomorphic functional calculus. If $a$ is an element of a Banach algebra, and $f$ is a holomorphic function on some neighborhood of the spectrum of $a$ we define
$$f(a):=\int_\gamma\frac{f(z)}{z-a}dz$$
where $\gamma$ is some path in the domain of $f$ that encircles the spectrum of $a$ and the integral is defined as the Riemann integral.
\end{definition}

The following basic facts follow easily from the definition or the cited Theorems.

\begin{proposition}\label{basicHolo}

\begin{enumerate}
    \item The map $f\mapsto f(a)$ is a morphism of Banach algebras. \cite[Theorem 3.3.5]{funementals1}
    \item If $a$ is a normal element of a $C^*$-algebra, then holomorphic functional calculus corresponds to the continuous functional calculus. (Follows from \cite[Theorem 3.3.5, Proposition 3.3.9]{funementals1} and Runge's Theorem.)
    \item If $\varphi$ is a morphism of Banach algebras, $\varphi(f(a))=f(\varphi(a))$.
    \end{enumerate}
\end{proposition}

\begin{definition}
A fiber bundle is said to be {\em numerable} if it is trivialized by a locally finite open cover that admits a partition of unity.
\end{definition}

\begin{definition}
A {\em universal principal $\Gamma$-bundle} is a numerable principal $\Gamma$-bundle $\alpha$ over a space $B\Gamma$ if
\begin{enumerate}
    \item if $\omega$ is another numerable principal $\Gamma$-bundle over a space $X$ there is a a continuous function $f:X\rightarrow B\Gamma$ so that $f^*(\alpha)$ is isomorphic to $\omega$.
    \item If $f$ and $g$ are both continuous functions from a space $X$ to $B\Gamma$ with so that $f^*(\alpha)$ is isomporphic to $g^*(\alpha)$ then $f$ is homotopy equivalent to $g$.
\end{enumerate}
For a group $\Gamma$ we define the {\em classifying space}, $B\Gamma$, to be the base space of a universal bundle; the total space will be denoted as $E\Gamma$.
\end{definition}

One can see immediately from the definition that $B\Gamma$ is unique up to homotopy. There is at least one model of the classifying space that is a CW-complex~\cite[Example 1b.7]{hatcher}.

The Mishchenko line bundle is the $\ell^1(\Gamma)$-bundle over $B\Gamma$ defined by $\ell_{B\Gamma}=E\Gamma\times_\Gamma \ell^1(\Gamma)$. More specifically, the total space $E\Gamma\times \ell^1(\Gamma)$ quotiented out by the diagonal action of $\Gamma$. The quotient map to $B\Gamma$ is given by $(a,b)\mapsto\pi(a)$ where $\pi$ is the map from $E\Gamma$ to $B\Gamma$.

\begin{lemma}\label{BGs}
Suppose that $X$ and $X'$ are both models of $B\Gamma$. Then there is a homotopy equivalence $h:X\rightarrow X'$ so that $\ell_{X}=h^*(\ell_{X'})$.
\end{lemma}
\begin{proof}
Let $Y$ and $Y'$ be models of $E\Gamma$, making universal bundles over $X$ and $X'$, respectively. By definition, there is a commutative diagram as follows:
\begin{center}
\begin{tikzcd}
Y\arrow{r}{h'}\arrow{d}{}&Y'\arrow{d}{}\\
X\arrow{r}{h}&X'
\end{tikzcd}
\end{center}
where $h'$ commutes with the action of $\Gamma$, and $h$ is a homotopy equivalence. It follows that we may make a diagram as follows
\begin{center}
\begin{tikzcd}
Y\times_{\Gamma}\ell^1(\Gamma)\arrow{r}{}\arrow{d}{}&Y'\times_\Gamma \ell^1(\Gamma)\arrow{d}{}\\
X\arrow{r}{h}&X'
\end{tikzcd}
\end{center}
One then sees by definition that the bundles on the left and right are $\ell_X$ and $\ell_{X'}$ respectively and that $h^*(\ell_{X'})=\ell_X$.
\end{proof}

Below we give $\ell^1(\Gamma)$ the structure of a Banach algebra where the multiplication is convolution.

\begin{definition}\label{q}
If $Y$ is a topological space with a continuous map $f:Y\rightarrow B\Gamma$ we may define $\ell_Y=f^*(\ell_{B\Gamma})$. Suppose that $\{U_1\ldots, U_N\}$ is a finite open cover of $Y$ so that $\ell_Y|U_i$ is trivial for each $i$. In this case, we represent $\ell_Y$ as a projection:
$$q=\sum_{i,j}\chi_i\chi_j\tensor g_{ij}\tensor e_{ij}\in C_b(Y)\tensor \ell^1(\Gamma)\tensor M_N.$$
where $\chi_i$ is a partition of unity in the sense that $\sum_{i=1}^N\chi_i^2=1$ and $g_{ij}\in\Gamma$ is a \v Cech cocycle representing $\ell_Y$.
\end{definition}

\begin{proposition}
If $q$ is defined as in Definition~\ref{q} then the right $C_b(Y,\ell^1(\Gamma))$-module of bounded sections on $\ell_Y$ is isomorphic to $pC_b(Y,\ell^1(\Gamma))^N$.
\end{proposition}

\begin{proof}
We may assume without loss of generality that the support of $\chi_i$ is $U_i$. If not, we note that $\{\supp(\chi_i)\}$ is also a finite open cover that trivializes $\ell_Y$, and leads to an identical definition of $q$.

One may identify sections of $\ell_Y$ with families of continuous functions $s_i:U_i\rightarrow \ell^1(\Gamma)$ that satisfy $s_i=g_{ij}s_j$. One checks that this condition is equivalent to
$$s_i=\sum_{j=1}^N\chi_j^2g_{ij} s_j.$$
One checks again that the above is equivalent to the condition
$$q\sum_j\chi_j s_j e_j=\sum_i\chi_i s_ie_i.$$
Thus $\{s_i\}\mapsto\sum_i\chi_i s_i e_i$ is the desired isomorphism.
\end{proof}

\begin{definition}\label{E_rho}
Let $\rho_n:\Gamma\rightarrow U(m_n)$ be asymptotically multiplicative in the operator norm. Let $f: Y\rightarrow B\Gamma$ be a continuous map and suppose that $\{U_i\}$, $\{\chi_i\}$ and $q$ are defined as in Definition~\ref{q}. Then $\rho_n$ can be extended to be a continuous linear function from $\ell^1(\Gamma)$ to $M_{m_n}$. This allows us to define $\rho_n(q)\in C_b(Y,M_{m_n})$. Furthermore we us holomorphic functional calculus to define the idempotent $\rho_{n\#}(q):=\chi_{\{z\in\C:\re(z)>1/2\}}(\rho_n(q))\in C_b(Y,M_{m_n})$, and we define $E_{\rho_n}^Y$ to be the vector bundle over $Y$ corresponding to this projection.
\end{definition}

A priori, $E_{\rho_n}^Y$ depends on our choice of partition of unity, but the following lemma shows that if $Y$ is compact, any two choices will eventually agree for a sufficiently multiplicative map.

\begin{lemma}\label{sameKT}
Let $Y$ be a compact Hausdorff space. Let $p$ and $q$ be any two idempotents in matrices over $C(Y,\ell^1(\Gamma))$ so that $pC(Y,\ell^1(\Gamma))^N$ and $qC(Y,\ell^1(\Gamma))^M$ are both isomorphic as right $C(Y,\ell^1(\Gamma))$-modules to sections on $\ell_Y$. For large enough $n$ that the bundles corresponding to $\rho_{n\#}(p)$ and $\rho_{n\#}(q)$ are isomporphic.
\end{lemma}

\begin{proof}
By \cite[Lemma~1.2.1]{rosenbergkt} there is some number $D\ge M,N$ so that if $p$ and $q$ are considered as elements of $M_D(C(Y,\ell^1(\Gamma)))$ then $p$ is conjugate to $q$ by an invertible element of $M_D(C(Y,\ell^1(\Gamma)))$. From this and \cite[Proposition 4.4.1]{blackadar} it follows that there is a homotopy of idempotents $p_t\in M_{2D}(C(Y,\ell^1(\Gamma)))$ connecting $p$ to $q$. For any $\varepsilon>0$ and any $(t,x)\in [0,1]\times Y$ there is some element $\gamma$ of $M_{2D}(\C[\Gamma])$ so that $||p_t(x)-\gamma||_1<\varepsilon$. By compactness of $[0,1]\times Y$ we may use finitely many elements of $\gamma_1,\ldots,\gamma_d$ to approximate $p_t(x)$ for all $(t,x)\in [0,1]\times Y$ up to $\varepsilon$. Then if $\rho_n$ is sufficiently multiplicative on the support of each $\gamma_i$ it follows that $||\rho_n(p_t)^2-\rho_n(p_t)||_\infty<\frac13$. Using holomorphic functional calculus, we can find a homotopy of idempotents from $\rho_{n\#}(p)$ to $\rho_{n\#}(q)$.
\end{proof}

\begin{proposition}\label{functorial}
Suppose that $f:Y\rightarrow Z$ and $g:Z\rightarrow B\Gamma$, and both $Y$ and $Z$ are compact. Suppose that $\rho_n:\Gamma\rightarrow U(m_n)$ is asymptotically multiplicative in the operator norm. Then for large enough $n$ we have $E_{\rho_n}^Y\cong f^*(E_{\rho_n}^Z)$.
\end{proposition}

\begin{proof}
Define $q_Y$ and $q_Z$ to represent $\ell_Y$ and $\ell_Z$ as in Definition~\ref{q}. Then by Remark~\ref{sameKT} we have that for large enough $n$, $\rho_{n\#}(q_Y)$ generates the same bundle as $\rho_{n\#}(f^*(q_Z))$. By functoriality of holomorphic functional calculus $\rho_{n\#}(f^*(q_Z))=f^*(\rho_{n\#}(q_z))$.
\end{proof}

\begin{remark}
Suppose that $\rho_n$ and $\psi_n$ are two operator norm asymptotic representations so that for all $g\in\Gamma$, $||\rho_n(g)-\psi_n(g)||_\infty\rightarrow0$. Then it follows that $||\rho_{n\#}(q)-\psi_{n\#}(q)||_\infty\rightarrow0$, which means that for large enough $n$ they represent the same class in K-theory. 
\end{remark}

We will use a formula for the Chern character that comes from Chern-Weil theory, and we will briefly explain some basics here. We will use the same notation as Appendix C in \cite{milnorstasheff}.

If $E$ is a smooth complex vector bundle, with fiber $V$, over a manifold $M$, and let $C^\infty(E)$ denote the space of smooth of $M$ sections of $E$. A connection is a map $\nabla$ from $C^\infty(E)$ to $\Lambda^1(M)\tensor C^\infty(E)$ that satisfies $\nabla(fs)=df\tensor s+f\nabla(s)$ for all $f\in C^\infty(M)$ and $s\in C^\infty(E)$. This induces a map $\hat{\nabla}:\Lambda^1(M)\tensor C^\infty(E)\rightarrow\Lambda^2(M)\tensor C^\infty(E)$ by the formula $\hat{\nabla}(\omega\tensor s)=d\omega\tensor s-\omega\wedge\nabla(s)$. Then the curvature associated with the $\nabla$ is the composition $K_\nabla=\hat{\nabla}\circ\nabla$. At each point, $x$, one may view $K_\nabla$ as a function from $V$ to $\Lambda^2(M,x)\tensor V$; this is because $K_\nabla$ evaluated at any point depends only on the value of the section at that point, not any derivatives or nearby points~\cite[Lemma 5]{milnorstasheff}. Thus, at each point, we may express $K_\nabla$ as a matrix of 2-forms that is unique up to a change of basis of $V$. This allows us to apply any basis-independent polynomial to $K_\nabla$, such as the trace or determinant.

It is well-known that the $k$th Chern character of a vector bundle can be written as $\frac{1}{(2\pi i)^kk!}\Tr(K_\nabla^k)$ in the $2k$ de Rham cohomology of $M$ \cite{milnorstasheff}.

Suppose that the vector bundle can be written as $pC(M)^n$ for some idempotent $p\in M_n(C^\infty(M))$. It is well-known and easy to check that $\nabla(s)=p(ds)$ is a connection. One checks that $\hat{\nabla}(dx\tensor s)=p d((dx)\cdot s)$. Using the product rule on $ds=d(ps)$, one can use the above to check that $K_\nabla(s)=pd(pd(ps))=(pdpdp)s$. To summarize the above discussion, we state the following.

\begin{proposition}\label{chernFormula}
Suppose that $E$ is a smooth vector bundle over $M$ so that sections of $E$ can be written as $pC_b(M)^n$ for some projection $p\in M_n(C_b^\infty(M))$. Then $\ch_k(M)$ is represented by the 2$k$-form $\frac{1}{(2\pi i)^kk!}\Tr((pdpdp)^k)$ the de Rham cohomology of $M$.
\end{proposition}

\section{Matrices of Differential Forms}

Let $M$ be a manifold. For $x\in M$ let $\Lambda^k(M,x)$ denote the $k$th exterior power of the complexification of the cotangent space at $x$. If $g$ is a Riemannian metric, $g$ induces a Hilbert space structure on the tangent space at each $x$, which induces a Hilbert space structure on the cotangent space, which induces a Hilbert space structure on~$\Lambda^k(M,x)$. Let $\Lambda_b^k(M)$ denote the space of continuous sections on the $k$th exterior power of the complexification of the cotangent bundle, that are bounded with respect to the norm induced by $g$.
\begin{definition}\label{norms}
Let $a\in M_n(\Lambda^k(M,x))$. For each $a$ there is a corresponding endomorphism $\Tilde{a}$ of the vector space $\bigoplus_{j=1}^{\dim(M)}\Lambda^j(M,x)^n$, defined by left multiplication by $a$. Using the Hilbert space structure induced by $g$, we may thus define the operator norm 
$$||a||_\infty=||\Tilde{a}||_\infty$$
or unnormalized Schatten $p$-norm 
$$||a||_p=\frac{||\Tilde{a}||_p}{2^{\dim(M)/p}}.$$ 
Note that these norms do depend on the choice of the Riemannian metric. For $a\in M_n(\Lambda^k(M))$ and $1\le p\le\infty$ we define
$$||a||_p=\sup_{x\in M}||a(x)||_p.$$
\end{definition}

\begin{remark}
Note that if $a\in M_n(C_b(M))\cong M_n(\Lambda_b^0(M))$ we note that the norm $||a||_p$ from Definition~\ref{norms} is equivalent to $\sup_{x\in M}||a(x)||_p$ since
$$\Tilde{a}(x)=a(x)\tensor\id_V$$
where $V=\bigoplus_{k=0}^{\dim(M)}\Lambda^k(M,x)$.
\end{remark}

\begin{theorem}[H\"older's Inequality]
Let $a\in M_n(\Lambda_b^k(M))$, $b\in M_n(\Lambda_b^j(M))$, $1\le p,q<\infty$ and $\frac1r=\frac1p+\frac1q$. The following hold
\begin{enumerate}
    \item $||ab||_r\le||a||_p||b||_q$
    \item $||ab||_p\le||a||_p||b||_\infty$
    \item $||ab||_p\le||a||_\infty||b||_p$.
\end{enumerate}
\end{theorem}

\begin{proof}
For matrices, part (1) follows from \cite[Theorem 2.8]{traceideas}. The other parts may be deduced from taking the limit as $q\rightarrow\infty$ or $p\rightarrow\infty$. One checks that for all $x\in M$ viewing $a(x)$ and $b(x)$ as endomorphisms of $\bigoplus_{j=1}^{\dim(M)}\Lambda^j(M,x)^n$ that composition of endomorphisms is the same as the matrix multiplication of $k$-forms. Finally one checks that multiplying both sides by $2^{-\dim(M)/p}$ or $2^{-\dim(M)/r}=2^{-\dim(M)/p-\dim(M)/q}$ does not change the inequalities.
\end{proof}

\begin{lemma}\label{tracebound}
Let $a\in M_n(\Lambda_b^k(M))$. Note that $\Tr(a)\in\Lambda_b^k(M)$, and $\Tr(a)(x)$ inherits a norm from Hilbert space structure induced by the Reimannian metric. Then $\sup_{x\in M}||\Tr(a)(x)||\le 2^{\dim(M)}\sqrt{\dim(M)\choose k}||a||_1$. 
\end{lemma}
\begin{proof}
Let $x\in M$. Let $N={\dim(M)\choose k}$ Then pick an orthonormal basis $\{v_i\}_{i=1}^N$ for $\Lambda^k(M,x)$. Let $\Tilde{a}(x)$ be the endomorphism of $\bigoplus_{j=0}^{\dim(M)}\Lambda^j(M,x)^n$ induced by multiplication by $a(x)$. Let $b$ be the restriction of $\tilde a(x)$ to $\Lambda^0(M,x)^n$ and let $\hat{v}_i$ be the linear function from $\C$ to $\Lambda^k(M,x)$ that maps 1 to $v_i$. There are $n\times n$ matrices $b_i$ so that we may write
$$b=\sum_{i=1}^N\hat{v}_i\tensor b_i.$$
Then we compute
$$b^*b=\sum_{i,j=1}^N\hat{v}_j^*\hat v_i\tensor b_j^*b_i=1_\C\tensor\sum_{i=1}^Nb_i^*b_i.$$
Because $b_i^*b_i$ are all positive, it follows that $1_\C\tensor b_i^*b_i\le b^*b$ and by~\cite[Theorem 2.2.6]{murph}, we may take the square root of both sides of this inequality. Thus
$$||b_i||_1=\Tr((b_i^*b_i)^{1/2})\le\Tr\left(\left(\sum_{i=1}^nb_i^*b_i\right)^{1/2}\right)=||b||_1.$$
Using the fact that $|\Tr(x)|\le||x||_1$ for any matrix $x$~\cite[Theorem 2.4.16]{murph}, we compute
\begin{align*}
||\Tr(a(x))||&=\left|\left|\sum_{i=1}^N\Tr(b_i)v_i\right|\right|\\
&=\sqrt{\sum_{i=1}^N|\Tr(b_i)|^2}\\
&\le\sqrt{\sum_{i=1}^N||b_i||_1^2}\\
&\le\sqrt{\sum_{i=1}^N||b||_1^2}\\
&=\sqrt{N}||b||_1\\
&\le\sqrt{N}||\Tilde{a}(x)||_1\\
&=2^{\dim(M)}\sqrt{N}||a(x)||_1.
\end{align*}
\end{proof}

\begin{definition}
We define $C_b^1(M)$ to be the space of once-differentiable bounded functions from $M$ to $\C$ whose exterior derivative is bounded in the norm induced by the Reimannian metric.
\end{definition}

\begin{definition}
We define the exterior derivative $d:M_n(C_b^1(M))\rightarrow M_n(\Lambda_{b}^{1}(M))$ by applying the (scalar) exterior derivative to each entry of the matrix.
\end{definition}

\begin{remark}
One can easily check from the definition of matrix multiplication that $d$ satisfies the following version of the product rule $d(ab)=d(a)b+ad(b)$.
\end{remark}

\begin{definition}\label{C1norm}
For $a\in M_n(C^1_b(M))$ define
$$||a||_{C^1}:=||a||_\infty+||da||_\infty.$$
\end{definition}

The following two theorems are likely to be well-known, but we provide proofs for the sake of completion.

\begin{theorem}\label{C^1banach}
The space $M_n(C_b^1(M))$ of matrices of $C^1$ functions from $M$ to $\C$ is a Banach algebra with respect to the norm defined in Definition~\ref{C1norm}.
\end{theorem}
\begin{proof}
To show this we note that $M_n(C_b^1(M))$ can be isometrically embedded in
$$M_n(C_b(M))\oplus M_n(\Lambda_b^1(M))$$
where the norm is the sum of the operator norm in each coordinate. The first coordinate is clearly complete. It is well known that continuous bounded sections of a locally trivial finite-dimensional vector bundle form a Banach space with the sup norm; one can prove this by embedding the space of sections in an $\ell^\infty$-product of spaces of bounded sections in trivial bundles. Thus $M_n(\Lambda_b^1(M))$ is a Banach space as well. The image of $M_n(C^1(M))$ under the embedding mentioned above is pairs $(a,b)$ so that $da=b$. We claim that this is true if and only if
\begin{equation}\label{integral}
\int_\gamma b=a(\gamma(\beta))-a(\gamma(\alpha))
\end{equation}
for all smooth paths $\gamma:[\alpha,\beta]\rightarrow M$. If $da=b$ then equation~\eqref{integral} follows from Stokes' theorem. Suppose equation~\eqref{integral} holds for all smooth paths. Let $x\in M$ We may pick a chart around $x$ and local coordinates $y_1,\ldots, y_{N}$. Without loss of generality, suppose that $x=0$ in this coordinate system. Let $\gamma(t)=e_it$. We use $\langle\cdot,\cdot\rangle$ to denote the formal pairing between cotangent vectors and tangent vectors.  We will compute the one-sided partial derivatives
\begin{align*}
\frac{\partial a_{jk}}{\partial y_i^+}|_x&=\lim_{h\rightarrow0^+}\frac{a_{jk}(\gamma(h))-a_{jk}(\gamma(0))}{h}\\
&=\lim_{h\rightarrow0^+}\frac1h\int_{\gamma|_{[0,h]}}b_{jk}\\
&=\lim_{h\rightarrow0^+}\frac1h\int_0^h\left\langle b_{jk}(\gamma(t)),\frac{\partial}{\partial y_i}\right\rangle dt\\
&=\left\langle b_{jk}(x),\frac{\partial}{\partial y_i}\right>.
\end{align*}
The computation for the partial derivatives coming from the negative direction is analogous. This shows that $\frac{\partial a}{\partial y_i}$ exists and that $da=b$ at each $x$. One can see that for each $\gamma$ the set of all pairs $(a,b)$ satisfying equation~\eqref{integral} is closed in the sup norm. Thus, as an intersection of closed sets, the set of points $(a,b)$ where $da=b$ is closed.

To show that $||\cdot||_{C^1}$ is submultiplicative, we compute
\begin{align*}
||ab||_{C^1}&=||ab||_\infty+||d(ab)||_\infty\\
&=||ab||_\infty+||d(a)b+abd||_\infty\\
&\le||a||_\infty||b||_\infty+||da||_\infty||b||_\infty+||a||_\infty||db||_\infty\\
&\le||a||_\infty||b||_\infty+||da||_\infty||b||_\infty+||a||_\infty||db||_\infty+||da||_\infty||db||_\infty\\
&= ||a||_{C^1}||b||_{C^1}.
\end{align*}
\end{proof}

\begin{theorem}\label{smooth}
Let $a\in M_n(\Lambda_b^0(M))$ be smooth. If $f$ is a holomorphic function on a neighborhood of the spectrum of $a$, then $f(a)$ is also smooth.
\end{theorem}

\begin{proof}
Suppose that $f$ is holomorphic on some open set $U\subseteq\C$. We will prove that the map $b\mapsto f(b)$ is a smooth map from $\{b\in M_n(\C):\sigma(b)\subseteq U\}$ to $M_n(\C)$. Note that
$$f(b)=\frac{1}{2\pi i}\int_\gamma\frac{f(z)}{(z-b)}dz.$$
Since this integral takes values in a finite-dimensional Banach space we may use the Lesbegue integral instead of the Riemann integral. Because inverting a matrix is smooth $f(b)$ is expressed as an integral of smooth functions. By a routine application of the dominated convergence theorem, we may pass any number of partial derivatives inside the integral.

Now let $a\in M_n(C^\infty(M))$. By assumption, the map $x\mapsto a(x)$ is smooth, and, as discussed above, $a(x)\mapsto f(a(x))$ is also smooth. Thus, $f(a)$ is a composition of smooth functions and is therefore smooth itself.
\end{proof}

\section{$p$-norm Almost Flatness}

\begin{definition}\label{notation}
Throughout this section, we will have the following notation. Fix $p\in\N_{>0}$. Let $\rho_n:\Gamma\rightarrow U(m_n)$ be a sequence of normalized functions that are asymptotically multiplicative in the unnormalized Schatten $p$-norm. Suppose that $\rho_n$ is normalized. Let $Y$ be a compact set with a smooth structure; in the proof of the main result, we will take $Y$ to be a compact subset of $\R^d$, but the proof also applies if $Y$ is a compact manifold with boundary. In the case that $Y\subseteq\R^n$ smooth maps are defined to be those that can be smoothly extended to a neighborhood of $Y$. Let $f: Y\rightarrow B\Gamma$ be continuous. Suppose that $M\subseteq Y$ is a manifold with a smooth structure given by the restriction of the smooth structure on $Y$, and suppose that the homology of $M$ is finitely generated. Let $\ell_Y$ be the pullback of the Mishchenko line bundle on $Y$.  Let $\{U_i\}_{i=1}^N$ be an open cover of $Y$, fine enough to observe local triviality of $\ell_Y$. Let $\{\chi_i\}_{i=1}^N$ be a smooth family of functions from $Y$ to $\R_{\ge0}$ that is the partition of unity in the sense that $\sum_{i=1}^N\chi_i^2=1$. Let $\{g_{ij}\}$ be a $\Gamma$ valued cocycle for $\ell_Y$ relative to the partition of unity $\{U_i\}_{i=1}^N$. Finally, let $e_{ij}$ be the $ij$th matrix unit. Following Definition~\ref{E_rho}, we define
$$\bar{a}_n:=\sum_{i,j=1}^N\chi_i\chi_j\tensor\rho_n(g_{ij})\tensor e_{ij}\in C^\infty(Y)\tensor M_{m_n}\tensor M_N.$$
We let $a_n$ be the restriction of $\bar{a}_n$ to $M$.
\end{definition}

\begin{lemma}\label{asquared}
$$\lim_{n\rightarrow\infty}||a_n^2-a_n||_p=0$$
and
$$\lim_{n\rightarrow\infty}||d(a_n^2-a_n)||_p=0.$$
\end{lemma}
\begin{proof}

Let
$$E_n:=\sum_{i,j,k=1}^N\chi_i\chi_j^2\chi_k\tensor(\rho_n(g_{ij})\rho_n(g_{jk})-\rho_n(g_{ik}))\tensor e_{ik}$$
We compute
\begin{align*}
a_n^2&=\sum_{i,j,\ell,k=1}^N\chi_i\chi_j\chi_\ell\chi_k\tensor\rho_n(g_{ij})\rho_n(g_{\ell k})\tensor e_{ij}e_{\ell k}\\
&=\sum_{i,j,k}\chi_i\chi_j^2\chi_k\tensor\rho_n(g_{ij})\rho_n(g_{jk})\tensor e_{ik}\\
&=\sum_{i,j,k}\chi_i\chi_j^2\chi_k\tensor\rho_n(g_{ik})\tensor e_{ik}+E_n\\
&=\sum_{i,k}\chi_i\chi_k\tensor\rho_n(g_{ik})\tensor e_{ik}+E_n\\
&=a_n+E_n
\end{align*}
So it suffices to show that the $p$-norms of $E_n$ and $d(E_n)$ converge to 0. This follows from the estimates
\begin{align*}
||E_n||_p&\le \sum_{i,j,k=1}^N\sup(\chi_i\chi_j^2\chi_k)\cdot||\rho_n(g_{ij})\rho_n(g_{jk})-\rho_n(g_{ik})||_p\\
||d(E_n)||_p&\le\sum_{i,j,k}^N\sup||d(\chi_i\chi_j^2\chi_k)||\cdot||\rho_n(g_{ij})\rho_n(g_{jk})-\rho_n(g_{ik})||_p.
\end{align*}
The suprema are finite since partial derivatives of $\chi_i$ can be extended to $Y$, which is compact.

\end{proof}

\begin{lemma}\label{bounded}
$$\limsup_{n\rightarrow\infty}||a_n||_\infty<\infty$$
and
$$\limsup_{n\rightarrow\infty}||da_n||_\infty<\infty.$$
\end{lemma}
\begin{proof}
For the first part, we estimate
$$||a_n||_\infty\le\sum_{i,j=1}^N\sup\chi_i\cdot\sup\chi_j\cdot||\rho_n(g_{ij})||_\infty\cdot||e_{ij}||_\infty\le N^2.$$
For the second part, we estimate
$$||d a_n||_\infty\le\sum_{i,j=1}^N\sup||d\chi_i||_\infty\sup||d\chi_j||_\infty$$
which is finite since partial derivatives of $\chi_i$ can be extended to $Y$, which is compact.
\end{proof}

\begin{definition}
Since $\rho_n$ is normalized, it is easy to compute that $a_n^*=a_n$. Since the operator norm is smaller than the $p$-norm, it follows that for large enough $n$ we may define a projection with continuous functional calculus as follows
$$q_n:=\chi_{[1/2,\infty)}(a_n).$$
By Proposition~\ref{basicHolo}, we may equivalently define $q_n$ by holomorphic functional calculus.
\end{definition}

Note that $q_n$ is smooth by Theorem~\ref{smooth}. This matches Definition~\ref{E_rho}.

\begin{lemma}\label{pnorm}
$$\lim_{n\rightarrow\infty}||a_n-q_n||_p=0.$$
\end{lemma}
\begin{proof}

Pick $y\in M$. Then we may choose an orthonormal basis of $\C^{m_n\cdot N}$ so that 
$$a_n(y)=\diag(\lambda_1,\ldots\lambda_{m_nN})$$
with respect to the chosen basis. Then $$q_n(y)=\diag(\chi_{[1/2,\infty)}(\lambda_1),\ldots,\chi_{[1/2,\infty)}(\lambda_{m_nN})).$$
Then using Lemma~\ref{calculation} we get
\begin{align*}
||a_n(y)-q_n(y)||_p&=||\langle\lambda_1-\chi_{[1/2,\infty)}(\lambda_1),\ldots,\chi_{[1/2,\infty)}(\lambda_{m_nN}\rangle||_p\\
&\le||\langle2(\lambda_1^2-\lambda_1),\ldots,2(\lambda_{m_nN}^2-\lambda_{m_nN})\rangle||_p\\
&=2||a_n(y)^2-a_n(y)||_p
\end{align*}
This goes to zero uniformly in $y$ by Lemma~\ref{asquared}.
\end{proof}
\begin{lemma}\label{operatornormD}
$$\lim_{n\rightarrow\infty}||d(a_n-q_n)||_\infty=0.$$
\end{lemma}
\begin{proof}
Let $E_n=a_n-q_n$, and $S=\{x\in\C:\re(x)\ge\frac12\}$. By the holomorphic functional calculus, we may express $E_n$ as the following integral:
$$E_n=\frac{1}{2\pi i}\int_{\gamma_r}\frac{z-\chi_S(z)}{z-a_n}dz$$
where $\gamma_r$ is a path parametrizing a pair of circles of radius $r$ around 0 and 1, and~$r$ is any number less than $\frac12$, but large enough that the spectrum of $a_n$ is encircled by $\gamma_r$. Note that the integral converges in $||\cdot||_{C^1}$ since $||\cdot||_{C^1}$ makes $C_b^1(M)$ a Banach algebra. Because $d$ is a bounded linear map with respect to $||\cdot||_{C^1}$ we may put $d$ inside the integral. One estimates that
\begin{align*}
||d(E_n)||&\le 2r\sup_{z\in\gamma_r}(|z-\chi_S(z)|\cdot||d(z-a_n)^{-1}||_\infty)\\
&=2r^2\sup_{z\in\gamma_r}||d(z-a_n)^{-1}||_\infty.
\end{align*}
We know that by Lemma~\ref{asquared} and the fact that $p$ norm dominates the operator norm, $||a_n^2-a_n||_{C^1}$ goes to zero. At the same time for $z\in\gamma_r$ we have $|z^2-z|=\min\{|z-1|,|z|\}\max\{|z-1|,|z|\}= r(1-r)$. So for large enough $n$, the sum
$$\sum_{j=0}^\infty(a_n^2-a_n)^j(z^2-z)^{-1-j}$$
converges, in $||\cdot||_{C^1}$, to an inverse to $z^2-z-a_n^2+a_n$. We note that $(z-a_n)(a_n+z-1)=z^2-z-a_n^2+a_n$, so it follows that
$$(z-a_n)^{-1}=(a_n+z-1)\sum_{j=0}^\infty(a_n^2-a_n)^j(z^2-z)^{-1-j}$$
One estimates that
\begin{align*}
||d(z-a_n)^{-1}||_\infty\le&||da_n||_\infty\sum_{j=0}^\infty||a_n^2-a_n||_\infty^j|z^2-z|^{-1-j}\\
&+||a_n+z-1||_\infty\sum_{j=1}^\infty j||d(a_n^2-a_n)||_\infty ||a_n^2-a_n||_\infty^{j-1}|z^2-z|^{-1-j}\\
\le&||da_n||_\infty(r(r-1))^{-1}\frac{1}{1-\frac{||a_n^2-a_n||_\infty}{r(1-r)}}\\
&+||a_n+z-1||_\infty||d(a_n^2-a_n)||_\infty\log\left(1-\frac{||a_n^2-a_n||_\infty}{r(r-1)}\right)
\end{align*}
By Lemma~\ref{bounded} and Lemma~\ref{asquared} this shows that there is a constant $C$, not depending on $r$ or $z$ so that
$$\limsup_{n\rightarrow\infty}||d(z-a_n)^{-1}||_\infty\le\frac{C}{r(r-1)}$$
Thus
$$\limsup_{n\rightarrow\infty}||d(E_n)||_\infty\le 2r^2\frac{C}{r(r-1)}=\frac{2rC}{r-1}.$$
By Lemma~\ref{calculation} and the spectral mapping theorem, for any given $n$, we may pick~$r$ to be anything less than $\frac12$ and greater than $2||a_n^2-a_n||_\infty$. By Lemma~\ref{asquared} and the fact that the $p$-norm dominates the $\infty$-norm, we have that we may pick $r$ to be constant greater than zero, and the above calculation will apply for sufficiently large $n$. Thus the desired result follows.
\end{proof}

\begin{lemma}\label{pnormD}
$$\lim_{n\rightarrow\infty}||d(a_n-q_n)||_p=0.$$
\end{lemma}
\begin{proof}
Let $a_n=q_n+E_n$. Note that $q_n$ and $E_n$ are both expressed in terms of the continuous functional calculus of $a_n$, so they commute with each other. We compute
\begin{align*}
d(a_n^2-a_n)&=d((q_n+E_n)^2-q_n-E_n)\\
&=d(q_n^2+2q_nE_n+E_n^2-q_n-E_n)\\
&=d(2q_nE_n+E_n^2-E_n)\\
&=d((2q_n-1)E_n)+d(E_n^2)\\
d(a_n^2-a_n)&=d(2q_n-1)E_n+(2q_n-1)dE_n+d(E_n)E_n+E_ndE_n
\end{align*}
We will now show that $||(2q_n-1)dE_n||_p\rightarrow0$ by arguing that each other term on both sides of the equation above goes to zero in the $p$-norm. The left side of the equation goes to zero by Lemma~\ref{asquared}. Moreover
$$||d(2q_n-1)E_n||_p\le2||dq_n||_\infty||E_n||_p.$$
By Lemma~\ref{asquared}, $||E_n||_p$ goes to zero, so it suffices to show that $||dq_n||_\infty$ is bounded. By Lemma~\ref{operatornormD} $dq_n$ is close in operator norm to $da_n$. By Lemma~\ref{bounded}, one sees that $||da_n||_\infty$ is bounded. Moreover
$$||d(E_n)E_n||_p\le ||dE_n||_\infty||E_n||_p$$
which goes to zero by Lemma~\ref{operatornormD} and Lemma~\ref{asquared}. The same argument applies to $E_ndE_n$. From this it follows that $||(2q_n-1)dE_n||_p$ goes to zero. Since $2q_n-1$ is a unitary, this equals $||dE_n||_p$ which thus goes to zero as well.
\end{proof}
\begin{lemma}\label{adada}
$$\lim_{n\rightarrow\infty}||a_nda_nda_n||_p=0.$$
\end{lemma}
\begin{proof}
Let
$$E_n:=\sum_{i,j,k=1}^N\chi_i(\chi_j^2d\chi_k+\chi_j\chi_kd\chi_j)\tensor(\rho_n(g_{ij})\rho_n(g_{jk})-\rho_n(g_{ik}))\tensor e_{ij}$$
Next, let
$$F_n:=\sum_{i,j,k,\ell=1}^N\chi_i(\chi_j^2\chi_kd\chi_kd\chi_\ell+\chi_j\chi_k^2d\chi_jd\chi_\ell+\chi_j\chi_k\chi_\ell d\chi_jd\chi_k)\tensor(\rho_n(g_{ik})\rho(g_{k\ell})-\rho_n(g_{i\ell}))\tensor e_{i\ell}$$
First, we see that
\begin{align*}
da_n=&\sum_{j,k=1}^N(\chi_jd\chi_k+\chi_kd\chi_j)\tensor\rho_n(g_{jk})\tensor e_{jk}\\
a_nda_n=&\sum_{i,j,k=1}^N\chi_i(\chi_j^2d\chi_k+\chi_j\chi_kd\chi_j)\tensor\rho_n(g_{ik})\tensor e_{ik}+E_n\\
a_nda_nda_n=&\sum_{i,j,k,\ell=1}^N\chi_i(\chi_j^2\chi_kd\chi_kd\chi_\ell+\chi_j\chi_k^2d\chi_jd\chi_\ell+\chi_j\chi_k\chi_\ell d\chi_jd\chi_k)\tensor\rho_n(g_{i\ell})\tensor e_{i\ell}\\
&+F_n+E_nda_n.
\end{align*}
We will show that the sum is zero by fixing $i$ and $\ell$. First, we see that
$$\sum_{j,k=1}^N\chi_j\chi_k d\chi_jd\chi_k=\sum_{1\le j<k\le N}(\chi_j\chi_k(d\chi_jd\chi_k+d\chi_kd\chi_j))=0.$$
Next note that
$$\sum_{j,k=1}^N(\chi_j^2\chi_kd\chi_k+\chi_j\chi_k^2d\chi_j)=\frac12\sum_{j,k=1}^Nd(\chi_j^2\chi_k^2)=\frac12d\left(\sum_{j,k=1}^N\chi_j^2\chi_k^2\right)=\frac12 d(1)=0.$$
We conclude that 
$$a_nda_nda_n=F_n+E_nda_n.$$
One estimates that
$$||E_n||_p\le\sum_{i,j,k=1}^N\sup||\chi_i(\chi_j^2d\chi_k+\chi_j\chi_kd\chi_j)||\cdot||\rho_n(g_{ij})\rho_n(g_{jk})-\rho_n(g_{ik}))||_p$$
and
$$||F_n||_p\le\sum_{i,j,k,\ell=1}^N\sup||\chi_i(\chi_j^2\chi_kd\chi_kd\chi_\ell+\chi_j\chi_k^2d\chi_jd\chi_\ell+\chi_j\chi_k\chi_\ell d\chi_jd\chi_k)||\cdot||\rho_n(g_{ik})\rho(g_{k\ell})-\rho_n(g_{i\ell})||_p$$
which both go to zero by asymptotic multiplicativity of $\rho_n$. Moreover, $||E_nda_n||_p\le||E_n||_p||da_n||_\infty$ which goes to zero by the above and Lemma~\ref{bounded}.
\end{proof}

\begin{theorem}\label{almostflat}
Let $\rho_n$, $\Gamma$, $Y$, $M$, and $q_n$ be as described in Definition~\ref{notation} The bundle $E_{\rho_n}^{M}$ is almost flat in the unnormalized Schatten $p$-norm in the sense that the curvature form, $q_ndq_ndq_n$, satisfies the following
$$\lim_{n\rightarrow\infty}||q_ndq_ndq_n||_p=0.$$
\end{theorem}
\begin{proof}
Let $q_n=a_n+E_n$. Then
$$q_ndq_ndq_n=E_ndq_ndq_n+a_ndE_ndq_n+a_nda_ndE_n+a_nda_nda_n$$
Thus
$$||q_ndq_ndq_n||_p\le||E_n||_p||dq_n||_\infty^2+||dE_n||_p||a_n||_\infty||dq_n||_\infty+||dE_n||_p||a_n||_\infty||da_n||_\infty+||a_nda_nda_n||_p$$
The first term on the right side goes to zero by Lemma~\ref{pnorm}, Lemma~\ref{bounded}, and Lemma~\ref{operatornormD}. The next two terms go to zero by Lemma~\ref{pnormD}, Lemma~\ref{bounded}, and Lemma~\ref{operatornormD}. The last term goes to zero by Lemma~\ref{adada}.
\end{proof}

\begin{theorem}\label{main2}
Let $\Gamma$, $p$, $\rho_n$, $Y$, and $M$ be as in Definition~\ref{notation}. Then for large enough $n$, $\ch_p(E_{\rho_n}^{M})=0$.
\end{theorem}

\begin{proof}
By Proposition~\ref{chernFormula}, $\ch_p(E_{\rho_n}^{M})$ is given by the de Rham cohomology class of $\omega_n=\Tr((q_ndq_ndq_n)^p)$. By H\"older's inequality, and Lemma~\ref{tracebound} we compute that
\begin{align*}
||\Tr((q_ndq_ndq_n)^p)||&\le2^{\dim(M)}\sqrt{\dim(M)\choose 2p}||(q_ndq_ndq_n)^p||_1\\
&\le2^{\dim(M)}\sqrt{\dim(M)\choose 2p}||q_ndq_ndq_n||_p^p
\end{align*}
which goes to zero by Theorem~\ref{almostflat}. If any singular $2p$ cohomology class $c$ may be expressed in terms of smooth maps $\Delta_j$ from the abstract $2p$-simplex to $Y$ as follows
$$c=\sum_{i=1}^Na_j\Delta_j\in Z_{2p}(M;\Z).$$
Note that $\Delta_j$ induces a linear map $\Delta_{j,x}^{2p}$ from $\Lambda^{2p}(Y,\Delta_j(x))$ to $\Lambda^{2p}(S_{2p},x)$ where $S_{2p}$ is the $2p$-simplex and $x\in S_{2p}$. Then, this map has an operator norm if both exterior spaces are equipped with Riemannian metrics. Since the operator norm depends continuously on $x$ it has some upper bound $C_j$. Then the pairing between de Rham cohomology and singular homology can be expressed as follows~\cite[ch 4.17]{deRhamThm}
\begin{align*}
|\langle\omega_n,c\rangle|&=\left|\sum_{j=1}^Na_j\int\Delta_j^*(\omega_n)\right|\\
&\le\sum_{j=1}^Na_jC_j||\omega_n||\mu(S_{2p})
\end{align*}
which goes to zero. However, $\langle\omega_n,c\rangle\in\frac{1}{p!}\Z$ since the $p$th Chern character of any vector bundle is in $\frac{1}{p!}$ times the image of $H^{2p}(M;\Z)$. Thus, for large enough $n$, $\langle\omega_n,c\rangle=0$. Since the homology of $M$ is finitely generated, there is some $n$ so that $\langle\omega_n,g\rangle=0$ for each generator of the $2p$ homology of $M$. By the universal coefficient theorem, we conclude that the real cohomology class of $\omega_n$ is zero.
\end{proof}

\section{Proof of Main Results}

We now prove Theorem~\ref{main},

\begin{proof}
Let $\rho_n:\Gamma\rightarrow U(m_n)$ be asymptotically multiplicative in the Schatten $p$-norm. By Proposition~\ref{balanced}, we may assume that $\rho_n$ is balanced. If $k\ge p$, then $\rho_n$ is also asymptotically multiplicative in the Schatten $k$-norm. First, pick $X$ to be a model of $B\Gamma$ that is a CW complex. Let $Y\subseteq X$ be compact. Since $Y$ is compact, it intersects only finitely many cells of $X$ \cite[Proposition A.1]{hatcher}. Thus, there is a finite CW complex $Y'$ satisfying $Y\subseteq Y'\subseteq B\Gamma$. For large enough $n$ we have that $E_{\rho_n}^{Y'}$ exists and $[E_{\rho_n}^Y]_0=\iota^*([(E_{\rho_n}^{Y'}]_0)$ by Proposition~\ref{functorial}. Since $Y'$ is a finite CW complex, there is an embedding $\iota: Y'\rightarrow\R^d$ so that there is an open neighborhood $Z$ of $\iota(Y')$ such that $\iota(Y')$ is a retract of $Z$ \cite[Corollary A.10]{hatcher}. By possibly picking a smaller $Z$ we may assume that $Z$ is a finite union of open balls and that $\iota(Y')$ is also a retract $\overline{Z}$. By the Mayer-Vietoris sequence, it follows that the homology of $\overline{Z}$ is finitely generated. By Letting $r: \overline{Z}\rightarrow Y'$ be the retraction we see that $[(E_{\rho_n}^{\overline{Z}})]_0=r^*([E_{\rho_n}^{Y'}]_0)$ by Proposition~\ref{functorial}. Let $\iota_Z$ be the inclusion of $Z$ to $\overline{Z}$. By Theorem~\ref{main2} we know that the restriction of $\iota_Z^*(\ch_k(E_{\rho_n}^{\overline{Z}}))=0$ for large enough $n$. By the naturality of Chern characters, and the fact that $\iota_Z^*\circ r^*$ is an injection of cohomology, it follows that $\ch_k(E_{\rho_n}^{Y'})=0$ as well. Finally, using the naturality of Chern characters again, we see that $\ch_k(E_{\rho_n}^Y)=0$. 

Now suppose that $X$ is a model for $B\Gamma$ that is not a CW complex and $Y\subseteq X$ is compact. We let $X'$ be a CW model for $B\Gamma$. There is homotopy equivalence $h:X\rightarrow X'$ so that $h^*(\ell_{X'})=\ell_X$ by Proposition~\ref{BGs}. If $Y'=h(Y)$, then one sees that $\ell_Y=h^*(\ell_{Y'})$. Then for large enough $n$, by Proposition~\ref{functorial} we see that $[E_{\rho_n}^Y]_0=[E_{\rho_n}^{Y'}]_0$. By the naturality of the Chern character and the above paragraph, we conclude that $\ch_k(E_{\rho_n}^Y)=0$ for sufficiently large $n$.
\end{proof}

\begin{remark}\label{chernclass}
Theorem~\ref{main} would not hold if Chern character were replaced with the Chern class. To see this, note that by~\cite[Proposition 3.3, Theorem 1.1]{partB} the first Chern class is given by the winding number invariant of Kazhdan, Exel, and Loring. Many examples of maps that are asymptotically multiplicative in the Schatten $p$-norm for all $p>1$ with nontrivial first Chern class can be found this way~\cite{KvoiculescuM}\cite{ELVoiculescuM}\cite{stab}\cite{frob}\cite{cup}. If two such examples exist so that $\ch_1(E_{\rho_n}^Y)\smile\ch_1(E_{\psi_n}^Y)\ne0$ one can use the Whitney sum formula to show that $c_2(E_{\rho_n\oplus\psi_n}^Y)=c_2(E_{\rho_n}^Y\oplus E_{\psi_n}^Y)\ne0$. One can easily make such examples by taking the direct product of any two groups mentioned in the above examples and applying the K\"uneth formula. The easiest to see explicit example is $\Z^4$. If $\rho_n$ takes the first two generators to the Voiculescu unitaries, and $\psi_n$ takes the remaining two generators to the Voiculescu unitaries. Then $c_2(E^{\mathbb{T}^4}_{\rho_n\oplus\psi_n})\ne0$, but $\rho_n\oplus\psi_n$ is asymptotically multiplicative in the Schatten $p$-norm for all $p>1$.
\end{remark}

Next we prove Corollary~\ref{dimbound}

\begin{proof}
Let $Y\subseteq B\Gamma$ be compact and let $M,k\in\N_{>0}$. Suppose, for contradiction, that the Corollary is false. Then let $\{S_m: m\in\mathbb {N} \}$ be an increasing family of finite subsets of $\Gamma$ whose union is all of $\Gamma$, and let $\varepsilon_m$ be a sequence of positive numbers tending towards 0. By the assumption that the Corollary is false, there is some $n_m\in\N$ and $\rho_m:\Gamma\rightarrow U(Mn_m^k)$ so that $\ch_k(E_{\rho}^Y)\ne0$ and for all $g,h\in S_m$, $||\rho_m(g)\rho_m(h)-\rho_m(gh)||_\infty<\frac{\varepsilon_m}{n_m}$. Then
$$||\rho_m(g)\rho_m(h)-\rho_m(gh)||_k\le M^{1/k}n_m||\rho_m(g)\rho_m(h)-\rho_m(gh)||_\infty<M^{1/k}\varepsilon_m.$$
This goes to zero; since all $g,h$ are eventually in $S_m$, this implies that $\rho_m$ is asymptotically multiplicative in the $k$-norm and so by Theorem~\ref{main} we get that $\ch_k(E^Y_{\rho_m})$ is eventually zero.
\end{proof}

The following is implicitly proved in the proof of the main result of~\cite{obs}. Since it is not explicitly stated there, we will include a proof that follows Dadarlat's argument there. 

\begin{theorem}[Dadarlat~\cite{obs}]\label{obsmain}
Let $\Gamma$ be a quasidiagonal group with a $\gamma$-element. If $H^{2k}(\Gamma;\Q)\ne0$ then there is a compact $Y\subseteq B\Gamma$ and an operator norm asymptotic homomorphisms $\varphi_n$ so that $\ch_k(E_{\rho_n}^Y)\ne0$.
\end{theorem}

First, we outline some notation; see~\cite{obs} for a full description. We let $\mathcal{Q}$ denote the universal UHF C*-algebra, $K(H)$ denote the compact operators on a separable infinite-dimensional Hilbert space, $M(K(H)\tensor\mathcal{Q})$ is the multiplier algebra, and $KK(C^*(\Gamma),\mathcal{Q})_{qd}$ is the subgroup of Kasparov's KK-theory consisting of Cuntz pairs $\varphi,\psi:C^*(\Gamma)\rightarrow M(K(H)\tensor\mathcal{Q})$ so that there is an increasing approximate unit of projections $p_n\in K(H)$ with $||[\psi(a),p_n\tensor 1_{\mathcal{Q}}]||\rightarrow0$ for all $a\in C^*(\Gamma)$.

\begin{proof}
By~\cite[Theorem 4.6]{obs}, the dual assembly map $\nu$ from $KK(C^*(\Gamma),\mathcal{Q})_{qd}$ to $\varprojlim(Y_i,\Q)$, where $Y_i$ is some sequence of finite subcomplexes who's union is $B\Gamma$, is surjective. Since the Chern character is a rational isomorphism~\cite{cherniso} for finite CW complexes, it follows that there is an element $x\in KK(C^*(\Gamma),\mathcal{Q})_{qd}$ and a compact $Y\subseteq B\Gamma$ so that $\ch_k(\nu_Y(x))\ne 0$. We may write $x$ as a Cuntz pair $(\varphi,\psi)$ where $\varphi$ and $\psi$ are nonzero $*$-homomorphisms from $C^*(\Gamma)$ to $M(K(H)\tensor\mathcal{Q})$ with the properties that $M(K(H)\tensor\mathcal{Q}))$ for all $a\in C^*(\Gamma)$, $\varphi(a)-\psi(a)\in K(H)\tensor\mathcal{Q}$ and there is an increasing approximate unit of projections $p_n\in K(H)$ so that $||[\varphi(a),p_n\tensor 1_{\mathcal{Q}}]||\rightarrow0$ for all $a\in C^*(\Gamma)$. It follows that  $||[\psi(a),p_n\tensor 1_{\mathcal{Q}}]||\rightarrow0$ as well. Let $\tilde\varphi_n=p_n\varphi p_n$ and $\tilde\psi_n=p_n\psi p_n$. We can find perturbations $\tilde\varphi_n$ and $\tilde\psi_n$ to ucp maps $\psi_n$ and $\varphi_n$ so that $\varphi_n(1)$ and $\psi_n(1)$ are projections as follows. We use functional calculus to find a projection $q_n$ close to $\tilde\varphi_n(1)$. Then $q_n\tilde\varphi_n(1)$ is invertible in $q_n(K(H)\tensor\mathcal{Q})q_n$ and define $\varphi_n=(q_n\tilde\varphi_n(1))^{-\frac12}\tilde\varphi_n(q_n\tilde\varphi_n(1))^{-\frac12}$. 

There are asymptotically multiplicative ucp maps $L_N:\mathcal{Q}\rightarrow M_{N!}$ defined by $L_N=\id_{M_1}\tensor\cdots\tensor\id_{M_N}\tensor\tau_{N+1}\tensor\tau_{N+2}\cdots$, where $\tau_N$ is the normalized trace on $M_N$. It follows from \cite[Proposition 2.5]{quasi-rep}, for large enough $n$,
$$\nu_Y(x)=[E_{\varphi_n}^Y]_0-[E_{\psi_n}^Y]_0\in K_0(C(Y)\tensor\mathcal{Q})\cong K^0(Y;\Q).$$
For large enough $N$ we will have that $\iota_*(L_{N\#}(\nu_Y(x)))=\nu_Y(x)$ where $\iota$ is the inclusion from $M_{N!}$ into $\mathcal{Q}$. Then for large enough $n,N$ we have
$$\ch_k(\nu_Y(x))=\ch_k(E_{L_N\circ\varphi_n}^Y)-\ch_k(E_{L_N\circ\psi_n}^Y).$$
In particular, at least one of the terms on the right is nonzero for all large enough $n$ and $N$; without loss of generality, suppose that $\ch_k(E_{L_N\circ\varphi_n}^Y)\ne0$. Picking some sequence $N_n$ so that $L_{N_n}\circ\varphi_n$ is asymptotically multiplicative, we have an asymptotic homomorphism with the desired properties.
\end{proof}
Corollary~\ref{can'timprove} follows from Theorem~\ref{main}, and Theorem~\ref{obsmain}.

Finally, we prove Corollary~\ref{cor3}. A key step is to cite the following result.

\begin{theorem}[\cite{partB} Theorem 1.2]\label{thm:unique}  Let $\Gamma$ be a torsion-free residually finite countable amenable group. For any finite set $F\subset \Gamma$ and any $\varepsilon>0,$ there exist a finite set $S\subset \Gamma$, $\delta>0$ and a compact subspace $Y\subset B\Gamma$
such that 
  for any two functions
$\rho,\rho':\Gamma \to U(k)$ with
$$||\rho(s)\rho(t)-\rho(st)||_\infty+||\rho'(s)\rho'(t)-\rho'(st)||_\infty<\delta\quad\forall s,t\in S$$
and $[E_\rho^Y]=[E_{\rho'}^Y]$ in $K^0(Y),$ there is a representation $\pi:\Gamma \to U(m)$ and a unitary $u\in U(k+m)$ such that 
\begin{equation*}
\|u(\rho(s)\oplus \pi(s))u^*-\rho'(s)\oplus \pi(s)\|_\infty<\varepsilon,\quad \forall s\in F.
\end{equation*}
\end{theorem}
\medskip

\noindent
{\em Proof of~Corollary~\ref{cor3}. } Let $\varepsilon>0$ and let $F\subseteq\Gamma$ be finite. Pick $\delta$, $S$ and $Y$ as in Theorem~\ref{thm:unique}. Eventually $\rho_n$ will satisfy
$$||\rho_n(s)\rho_n(t)-\rho_n(st)||_\infty<\delta\quad\forall s,t\in S$$
since the $p$-norm dominates the operator norm. Moreover, by Theorem~\ref{main}, for all $k>0$, $\ch_k([E_{\rho_n}^Y])=0$ for large enough $n$. Since the Chern character is a rational isomorphism~\cite{cherniso}, we may pick $r_n$ so that $r_n[E_{\rho_n}^Y]$ is isomorphic to a trivial bundle in $K^0(Y)$. Applying Theorem~\ref{thm:unique} to $\rho'=\rho_n^{\oplus r_n}$ and $\rho$ the $r_n\cdot m_n$th amplification of the trivial representation we get that $\rho'\oplus \pi$ is approximated up to $\varepsilon$ in operator norm on $F$. Since $\varepsilon$ can be made arbitrarily small and $F$ can be made arbitrarily large, the desired result follows.
\begin{flushright}
$\square$
\end{flushright}

\noindent
{\bf Acknowledgments:} I would like to thank Rufus Willett for many helpful conversations that pointed me in the right direction for several parts of this paper. I want to thank Thomas Hangelbroek for a conversation that helped develop the proof of Theorem~\ref{C^1banach}.

\bibliographystyle{plain}
\bibliography{main}
\end{document}